\documentclass[a4paper,fleqn,11pt]{article}
\usepackage{amsmath}
\usepackage{latexsym}
\usepackage{amssymb}
\usepackage{graphicx}
\usepackage{amsthm}

\newtheorem{thm}{Theorem}%[section]
\newtheorem{lem}[thm]{Lemma}
\newtheorem{cor}[thm]{Corollary}

\theoremstyle{definition}
\newtheorem{defn}{Definition}%[section]

\begin{document}

\title{A Linear-Time Algorithm for the Maximum Matched-Paired-Domination Problem in Cographs}

%\date{}
\author{\vspace{0.5cm}Ruo-Wei Hung\thanks{Corresponding author's e-mail: rwhung@cyut.edu.tw} and Chih-Chia Yao\\
\textit{Department of Computer Science and Information Engineering,}\\
\textit{Chaoyang University of Technology,}\\
\textit{Wufong, Taichung 41349, Taiwan}\\
E-mail addresses: \{rwhung, ccyao\}@cyut.edu.tw}

%%% ----------------------------------------------------------------------
\maketitle
%\thispagestyle{empty}
%%% ----------------------------------------------------------------------

\begin{abstract}
Let $G=(V,E)$ be a graph without isolated vertices. A matching in
$G$ is a set of independent edges in $G$. A perfect matching $M$
in $G$ is a matching such that every vertex of $G$ is incident to
an edge of $M$. A set $S\subseteq V$ is a
\textit{paired-dominating set} of $G$ if every vertex in $V-S$ is
adjacent to some vertex in $S$ and if the subgraph $G[S]$ induced
by $S$ contains at least one perfect matching. The
paired-domination problem is to find a paired-dominating set of
$G$ with minimum cardinality. In this paper, we introduce a
generalization of the paired-domination problem, namely the
maximum matched-paired-domination problem. A set $MPD\subseteq E$
is a \textit{matched-paired-dominating set} of $G$ if $MPD$ is a
perfect matching of $G[S]$ induced by a paired-dominating set $S$
of $G$. Note that the paired-domination problem can be regard as
finding a matched-paired-dominating set of $G$ with minimum
cardinality. Let $\mathcal{R}$ be a subset of $V$, $MPD$ be a
matched-paired-dominating set of $G$, and let $V(MPD)$ denote the
set of vertices being incident to edges of $MPD$. A
\textit{maximum matched-paired-dominating set} $MMPD$ of $G$
w.r.t. $\mathcal{R}$ is a matched-paired-dominating set such that
$|V(MMPD)\cap \mathcal{R}|\geqslant |V(MPD)\cap \mathcal{R}|$. An
edge in $MPD$ is called \textit{free-paired-edge} if neither of
its both vertices is in $\mathcal{R}$. Given a graph $G$ and a
subset $\mathcal{R}$ of vertices of $G$, the \textit{maximum
matched-paired-domination problem} is to find a maximum
matched-paired-dominating set of $G$ with the least
free-paired-edges; note that, if $\mathcal{R}$ is empty, the
stated problem coincides with the classical paired-domination
problem. In this paper, we present a linear-time algorithm to
solve the maximum matched-paired-domination problem
in cographs.\\


\end{abstract}

\newpage
\baselineskip=18.8pt
%====================================================================
\section{Introduction}\label{Introduction}
%====================================================================
All graphs considered in this paper are finite and undirected,
without loops or multiple edges. Let $G=(V,E)$ be a graph without
isolated vertices. The \textit{open neighborhood} $N_G(v)$ of the
vertex $v$ in $G$ is defined to be $N_G(v)=\{u\in V|uv\in E\}$ and
the \textit{closed neighborhood} $N_G[v]$ of $v$ is
$N_G(v)\cup\{v\}$. For a set $S\subseteq V$, the subgraph of $G$
induced by the vertices in $S$ is denoted by $G[S]$. A set
$D\subseteq V$ is a \textit{dominating set} of $G$ if every vertex
not in $D$ is adjacent to at least a vertex in $D$. The
\textit{domination problem} is to find a dominating set of $G$
with minimum cardinality. The bibliography in domination and its
variations maintained by Haynes et al. \cite{Haynes98a} currently
has over 1200 entries; Hedetniemi and Laskar \cite{Hedetniemi90}
edited a special issue of \textit{Discrete Mathematics} devoted
entirely to domination, and two books on domination and its
variations in graphs \cite{Haynes98a,Haynes98b} have been written.

\begin{figure}[t]
\begin{center}
\includegraphics[scale=0.7]{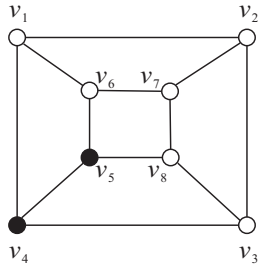}
\caption{The tree-cube graph $Q_3$, where restricted vertices are
drawn by filled circles.} \label{Fig_example_paired-domination}
\end{center}
\end{figure}

A \textit{matching} in a graph $G$ is a set of independent edges
in $G$. A \textit{perfect matching} $M$ in a graph $G$ is a
matching such that every vertex of $G$ is incident to an edge of
$M$. A \textit{paired-dominating set} of a graph $G$ is a set $PD$
of vertices of $G$ such that $PD$ is a dominating set of $G$ and
$G[PD]$ contains at least one perfect matching. In other words, a
paired-dominating set with matching $M$ is a dominating set
$PD=\{v_1, v_2, \cdots, v_{2t-1}, v_{2t}\}$ with independent edge
set $M=\{e_1, e_2, \cdots, e_t\}$, where each edge $e_i$ joins two
vertices of $PD$. The minimum cardinality of a paired-dominating
set for a graph $G$ is called the \textit{paired-domination
number}, denoted by $\gamma_{\rm p}(G)$. A paired-dominating set
of $G$ with cardinality $\gamma_{\rm p}(G)$ is called a
\textit{minimum paired-dominating set} of $G$. The \textit{paired
domination problem} is to find a minimum paired-dominating set of
$G$. Note that every graph without isolated vertices contains a
minimum paired-dominating set \cite{Haynes98c}. For example, for
the three-cube graph $Q_3$ in Fig.
\ref{Fig_example_paired-domination}, $PD=\{v_1, v_2 , v_3, v_4\}$
with matching $M_1=\{v_1v_2, v_3v_4\}$ or $PD$ with matching
$M_2=\{v_1v_4, v_2v_3\}$ is a minimum paired-dominating set of
$Q_3$ and $\gamma_{\rm p}(Q_3)=4$.

Paired-domination was introduced by Haynes and Slater and the
decision problem to determine $\gamma_{\rm p}(G)$ of an arbitrary
graph $G$ has been known to be NP-complete \cite{Haynes98c}. It is
still NP-complete on some special classes of graphs, including
bipartite graphs, chordal graphs, and split graphs \cite{Chen08}.
However, it admits polynomial time algorithms when the input is
restricted to be in some special classes of graphs, including
trees \cite{Qiao03}, circular-arc graphs \cite{Cheng07},
permutation graphs \cite{Cheng09}, block graphs, and interval
graphs \cite{Chen08}.

Paired-domination has found the following application
\cite{Haynes98c}. In a graph $G$ if we think of each vertex $s$ as
the possible location for a guard capable of protecting each
vertex in $N_G[s]$, then ``domination" requires every vertex to be
protected. In paired-domination, each guard is assigned another
adjacent one, and they are designed as backups for each other.
However, some locations may play more important role (for example,
important facilities are placed on these locations) and, hence,
they are placed by guards for instant monitoring and protection
possible. In this application, the number of guards placed on the
important locations is as large as possible. Motivated by the
above issue we introduce a generalization of the paired-domination
problem, namely, the maximum matched-paired-domination problem.

Let $G=(V, E)$ be a graph without isolated vertices, $\mathcal{R}$
be a subset of $V$, and let $PD$ be a paired-dominating set of
$G$. For a set $M$ of independent edges in $G$, we use $V(M)$ to
denote the set of vertices being incident to edges of $M$. A set
$MPD\subseteq E$ is called a \textit{matched-paired-dominating
set} of $G$ if $MPD$ is a perfect matching of $G[PD]$ induced by a
paired-dominating set $PD$ of $G$. That is, $V(MPD)$ is a
paired-dominating set $PD$ of $G$ and $MPD$ specifies a perfect
matching of $G[PD]$. Note that the paired-domination problem can
be regard as finding a matched-paired-dominating set of $G$ with
minimum cardinality. For an edge $e=uv\in MPD$, we say that $e$ is
a \textit{paired-edge} in $MPD$, $u$ is paired with $v$, and $u$
is the \textit{partner} of $v$. In addition, we will use $\langle
u, v \rangle$ to denote a paired-edge $uv$ in $MPD$ if it is
understood without ambiguity. Note that in a paired-dominating set
$PD$ of $G$, it is necessary to specify which vertex is the
partner of a vertex in $PD$. The \textit{matched number} of a
matched-paired-dominating set $MPD$ is defined to be $|V(MPD)\cap
\mathcal{R}|$. The \textit{maximum matched number} $\beta(G)$ of
$G$ is defined to be the largest matched number of a
matched-paired-dominating set in $G$. A \textit{maximum
matched-paired-dominating set} of $G$ w.r.t. $\mathcal{R}$ is a
matched-paired-dominating set with matched number $\beta(G)$. A
paired-edge in $MPD$ is called \textit{free-paired-edge} if both
of its vertices are not in $\mathcal{R}$. A
matched-paired-dominating set of $G$ is called \textit{canonical}
if it is a maximum matched-paired-dominating set of $G$ with the
least free-paired-edges. Given a graph $G$ and a subset
$\mathcal{R}$ of vertices of $G$, the \textit{maximum
matched-paired-domination problem} is to find a canonical
matched-paired-dominating set of $G$ w.r.t. $\mathcal{R}$. Note
that if $\mathcal{R}$ is empty, the stated problem coincides with
the classical paired-domination problem. We call $\mathcal{R}$ the
\emph{restricted vertex set} of $G$. The vertices in $\mathcal{R}$
are called \emph{restricted vertices} and the other vertices are
called \emph{free vertices}. For example, given a graph $G$ and a
restricted vertex set $\mathcal{R}=\{v_4, v_5\}$ shown in Fig.
\ref{Fig_example_paired-domination}, let $MPD_1=\{\langle v_1,v_2
\rangle, \langle v_3,v_4 \rangle\}$, $MPD_2=\{\langle v_4,v_5
\rangle, \langle v_2,v_7 \rangle\}$, and let $MPD_3=\{\langle
v_1,v_4 \rangle, \langle v_5,v_6 \rangle\}$. We can see that
$\beta(G)=2\leqslant |\mathcal{R}|$. Then, both $MPD_2$ and
$MPD_3$ are maximum matched-paired-dominating sets of $G$, but
$MPD_1$ is not a maximum matched-paired-dominating set of $G$. It
is straightforward to see that $MPD_2$ contains a free-paired-edge
and $MPD_3$ contains no free-paired-edge. Thus, $MPD_3$ is a
canonical matched-paired-dominating set of $G$, but $MPD_2$ is not
canonical.

Now, we review cographs. Cographs (also called
complement-reducible graphs) are defined as the class of graphs
formed from a single vertex under the closure of the operations of
\textit{union} and \textit{complement}. Cographs were introduced
by Lerchs \cite{Lerchs71}, who studied their structural and
algorithmic properties and enumerated the class. Names synonymous
with cographs include $D^*$-graphs, $P_4$ restricted graphs, and
Hereditary Dacey graphs. Several characterizations of cographs are
known. For example, it is shown that $G$ is a cograph if and only
if $G$ contains no $P_4$ (a path consisting of four vertices) as
an induced subgraph \cite{Corneil81}. Cographs have arisen in many
disparate areas of mathematics and have been independently
rediscovered by various researchers. These graphs can be
recognized in linear time \cite{Corneil85,Habib05}. The class of
cographs forms a subclass of distance-hereditary graphs
\cite{Corneil81,Corneil85} and permutation graphs, and is a
superclass of threshold graphs and complete-bipartite graphs.
Numerous properties and optimization problems in these graphs have
been studied
\cite{Bodlaender89,Bodlaender93,Chang08,Damiand01,Hsieh07,Hung06,
Larrion04,Nakano03,Nikolopoulos05,Retore03,Shamir04,Yu93}. In this
paper, we will solve the maximum matched-paired-domination problem
on cographs in linear time.

%=================================================================
\section{Known Results and Terminology}\label{Preliminaries}
%=================================================================
Let $G$ be a graph without isolated vertices. Haynes and Slater
showed that a paired-dominating set of $G$ does exist and
$\gamma_{\rm p}(G)$ is even \cite{Haynes98c}.

\begin{lem}\cite{Haynes98c}\label{paired-domination_exists}
Let $G$ be a graph without isolated vertices. Then, there exists a
paired-dominating set in $G$ and $\gamma_{\rm p}(G)$ is even.
\end{lem}
%\begin{proof}
%We prove the lemma by induction on the number of vertices of
%$G=(V,E)$. Initially, let $|V|=2$ and $V=\{v_1,v_2\}$. Since $G$
%has no isolated vertices, $v_1$ is adjacent to $v_2$. Let
%$PD=\{v_1, v_2\}$. Then, $PD$ is a paired-dominating set of $G$.
%Assume that for any graph $\widehat{G}=(\widehat{V},\widehat{E})$
%such that $|\widehat{V}|=k\geqslant 2$ and $\widehat{G}$ contains
%no isolated vertex, $\widehat{G}$ has a paired-dominating set. Let
%$G=(V,E)$ be a graph such that $|V|=k+1$ and $G$ contains no
%isolated vertex, and let $v\in V$. We denote by $G-v$ deleting $v$
%and edges incident to $v$ from $G$. Then, the number of vertices
%of $G-v$ equals $k$. By induction hypothesis, $G-v$ has a
%paired-dominating set $PD$. If $N_G(v)\subseteq PD$, then $PD$ is
%a paired-dominating set of $G$. Suppose $N_G(v)\not\subseteq PD$.
%Let $u$ be a vertex in $N_G(v)-PD$. Then, $PD \cup \{v, u\}$ is a
%paired-dominating set of $G$. Thus, $G$ has a paired-dominating
%set. By induction, a graph without isolated vertices has a
%paired-dominating set.
%\end{proof}

It follows from Lemma \ref{paired-domination_exists} that we have
the following corollary.

\begin{cor}\label{matched-paired-domination_exists}
Let $G$ be a graph without isolated vertices. Then, there exists a
canonical matched-paired-dominating set in $G$.
\end{cor}

The following lemma is easily verified from the definition.

\begin{lem}\label{trivial_maximum_matched-paired-domination}
Assume $G$ is a graph without isolated vertices and $\mathcal{R}$
is a restricted vertex set of $G$. Let $MPD$ be a
matched-paired-dominating set of $G$ w.r.t. $\mathcal{R}$. Then,\\
(1) if $|V(MPD)\cap \mathcal{R}|= |\mathcal{R}|$ and
$|MPD|=\lceil\frac{|\mathcal{R}|}{2}\rceil$, then
$\beta(G)=|\mathcal{R}|$ and $MPD$ is a canonical
matched-paired-dominating set of $G$;\\
(2) if $|V|=|\mathcal{R}|$ is odd, $|V(MPD)\cap \mathcal{R}|=
|\mathcal{R}|-1$, and
$|MPD|=\lfloor\frac{|\mathcal{R}|}{2}\rfloor$, then
$\beta(G)=|\mathcal{R}|-1$ and $MPD$ is a canonical
matched-paired-dominating set of $G$.
\end{lem}

Now, we define some notations to be used in the paper. In the
following, we use $\mathcal{R}$ to denote the restricted vertex
set of a graph $G$.

\begin{defn}
A paired-edge in a matched-paired-dominating set of $G$ w.r.t.
$\mathcal{R}$ is called \textit{full-paired-edge} if both of its
vertices are in $\mathcal{R}$, is called \textit{semi-paired-edge}
if its one vertex is in $\mathcal{R}$ but the other vertex is not
in $\mathcal{R}$, and is called \textit{free-paired-edge} if both
of its vertices are not in $\mathcal{R}$.
\end{defn}

\begin{defn}
A matched-paired-dominating set $MPD$ of $G$ w.r.t. $\mathcal{R}$
is called \textit{$(k,s,f)$-matched-paired-dominating set} if (1)
$|MPD|=k+s+f$; (2) there are exactly $k$ full-paired-edges in
$MPD$; (3) there are exactly $s$ semi-paired-edges in $MPD$, and
(4) all other paired-edges in $MPD$ are free-paired-edges.
\end{defn}

By the above definition, a paired-edge in a
$(k,s,f)$-matched-paired-dominating set $MPD$ is either a
full-paired-edge, a semi-paired-edge or a free-paired-edge. Then,
the matched number of $MPD$ is $|V(MPD)\cap\mathcal{R}|=2k+s$.
Thus, a maximum $(k^*,s^*,f^*)$-matched-paired-dominating set of a
graph $G$ satisfies that $\beta(G)=2k^*+s^*\geqslant 2k+s$ for any
$(k,s,f)$-matched-paired-dominating set of $G$.

\begin{defn}
Let $MPD$ be a $(k,s,f)$-matched-paired-dominating set of a graph
$G$ w.r.t. $\mathcal{R}$. Define $K_G(MPD)$, $S_G(MPD)$, and
$F_G(MPD)$ to be the subsets of $MPD$ consisting of all
full-paired-edges, all semi-paired-edges, and all
free-paired-edges in $MPD$, respectively.
\end{defn}

\begin{figure}[t]
\begin{center}
\includegraphics[scale=0.85]{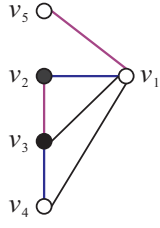}
\caption{A graph $G$ with restricted vertex set
$\mathcal{R}=\{v_2,v_3\}$, where restricted vertices are drawn by
filled circles.}\label{Fig_restricted_paired-domination}
\end{center}
\end{figure}

For example, let $G$ be a graph with restricted vertex set
$\mathcal{R}=\{v_2,v_3\}$ shown in Fig.
\ref{Fig_restricted_paired-domination}. Let $MPD_1=\{\langle
v_2,v_3 \rangle, \langle v_1,v_5 \rangle\}$ and let
$MPD_2=\{\langle v_1,v_2 \rangle, \langle v_3,v_4 \rangle\}$.
Then, $MPD_1$ is a $(1,0,1)$-matched-paired-dominating set and
$MPD_2$ is a $(0,2,0)$-matched-paired-dominating set. By
definition, $K_G(MPD_1)=\{\langle v_2,v_3 \rangle\}$,
$S_G(MPD_1)=\emptyset$, and $F_G(MPD_1)=\{\langle v_1,v_5
\rangle\}$, where $|K_G(MPD_1)|=1$, $|S_G(MPD_1)|=0$, and
$|F_G(MPD_1)|=1$.

Next, we introduce cographs. A graph is a cograph if there is no
induced path containing four vertices \cite{Corneil81}. Such
graphs are exactly the class of distance-hereditary graphs with
diameters less than or equal to two \cite{Bandelt86}. Every
cograph can be recursively defined as follows.

\begin{defn}\label{Cographs}\cite{Corneil81,Corneil85}
The class of cographs can be defined by the following recursive
definition:\\
(1) A graph consisting of a single vertex and no edges is a
cograph.\\
(2) If $G_L=(V_L,E_L)$ and $G_R=(V_R,E_R)$ are cographs, then the
\textit{union} $G$ of $G_L$ and $G_R$, denoted by $G=G_L\oplus
G_R=(V_L\cup V_R,E_L\cup E_R)$, is a cograph. In this case, we say
that $G$ is formed from $G_L$ and $G_R$ by a \textit{union
operation}.\\
(3) If $G_L=(V_L,E_L)$ and $G_R=(V_R,E_R)$ are cographs, then the
\textit{joint} $G$ of $G_L$ and $G_R$, denoted by $G=G_L\otimes
G_R=(V_L\cup V_R,E_L\cup E_R\cup \widehat{E})$, is a cograph,
where $\widehat{E}=\{uv|\forall u\in V_L$ and $\forall v\in
V_R\}$. In this case, we say that $G$ is formed from $G_L$ and
$G_R$ by a \textit{joint operation}.
\end{defn}

A cograph $G$ can be represented by a rooted binary tree $DT(G)$,
called a \emph{decomposition tree} \cite{Chang97,Corneil81}. The
leaf nodes of $DT(G)$ represent the vertices of $G$. Each internal
node of $DT(G)$ is labeled by either `$\oplus$' or `$\otimes$'.
The cograph corresponding to a $\oplus$-labeled (resp.
$\otimes$-labeled) node $v$ in $DT(G)$ is obtained from the
cographs corresponding to the children of $v$ in $DT(G)$ by means
of a \textit{union} (resp. \textit{joint}) operation. A
decomposition tree of a cograph can be constructed as follows.

\begin{defn}\cite{Chang97}\label{decomposition-tree}
The decomposition tree $DT(G)$ of a cograph $G$ consisting of a
single vertex $v$ is a tree of one node labeled by $v$. If $G$ is
formed from $G_L$ and $G_R$ by a union (resp. joint) operation,
then the root of the decomposition $DT(G)$ is a node labeled by
$\oplus$ (resp. $\otimes$) with the roots of $DT(G_L)$ and
$DT(G_R)$ being the children of the root of $DT(G)$, respectively.
\end{defn}

\begin{figure}[t]
\begin{center}
\includegraphics[width=0.7\textwidth]{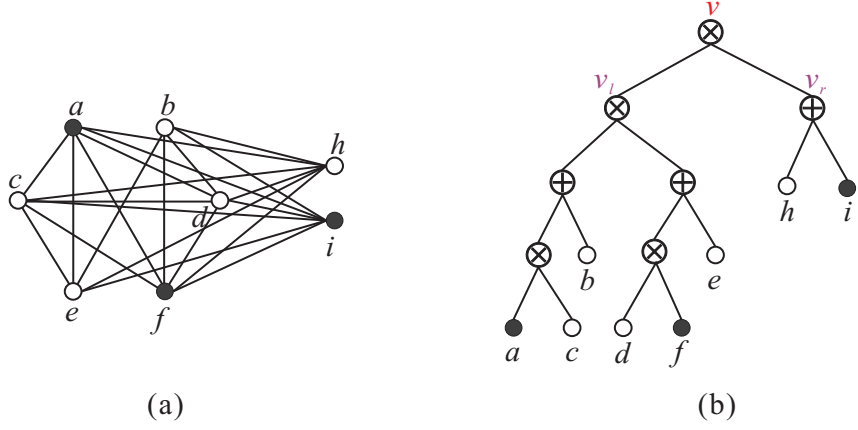}
\caption{(a) A cograph $G$ with restricted vertex set
$\mathcal{R}=\{a,f,i\}$, and (b) a decomposition tree $DT(G)$ of
$G$, where restricted vertices are drawn by filled
circles.}\label{Fig_decomposition-tree}
\end{center}
\end{figure}

The decomposition tree $DT(G)$ of a cograph $G$ is a rooted and
unordered binary tree. Note that exchanging the left and right
children of an internal node in $DT(G)$ will be also a
decomposition tree of $G$. For instance, given a cograph $G$ shown
in Fig. \ref{Fig_decomposition-tree}(a), the decomposition tree
$DT(G)$ of $G$ is shown in Fig. \ref{Fig_decomposition-tree}(b).

\begin{thm}\cite{Chang97,Corneil81}\label{decomposition-tree-theorem}
A decomposition tree $DT(G)$ of a cograph $G=(V,E)$ can be
constructed in $O(|V|+|E|)$-linear time.
\end{thm}

\section{The Maximum Matched-Paired-Domination Problem on\\ Cographs}

In this section, we will show that the maximum
matched-paired-domination problem on cographs is linear solvable.
Recall that a canonical matched-paired-dominating set of a graph
is a maximum matched-paired-dominating set with the least
free-paired-edges.  In fact, we will construct a canonical
matched-paired-dominating set of a connected cograph in linear
time. We first give the following lemma to show some properties of
a maximum matched-paired-dominating set of a graph.

\begin{lem}\label{Properties}
Assume $G$ is a connected graph without isolated vertices and
$\mathcal{R}$ is a restricted vertex set of $G$. Let $MMPD$ be a
maximum matched-paired-dominating set of $G$ w.r.t. $\mathcal{R}$
and let $\upsilon\in\mathcal{R}-V(MMPD)$. Then, the following
statements hold true:\\
(1) if $\langle v_f,\widetilde{v}_f \rangle$ is a free-paired-edge
in $MMPD$, then $\upsilon$ is adjacent to neither $v_f$ nor
$\widetilde{v}_f$;\\
(2) $N_G(\upsilon)\subseteq V(MMPD)$;\\
(3) if $\langle v_x,\widetilde{v}_x \rangle$ is a semi-paired-edge
or full-paired-edge in $MMPD$ and $\upsilon$ is adjacent to $v_x$,
then $N_G(\widetilde{v}_x)-\{\upsilon\}\subseteq V(MMPD)$;\\
(4) if $\langle v_f,v_r \rangle$ is a semi-paired-edge, with
restricted vertex $v_r$, in $MMPD$, then $\upsilon$ is not
adjacent to $v_r$.
\end{lem}
\begin{proof}
We first prove Statement (1). Assume by contradiction that
$\upsilon$ is adjacent to $v_f$. If
$N_G(\widetilde{v}_f)-V(MMPD)=\emptyset$, then $MMPD-\{\langle
v_f,\widetilde{v}_f \rangle\}\cup\{\langle v_f,\upsilon \rangle\}$
is a matched-paired-dominating set of $G$ having more restricted
vertices than $MMPD$, a contradiction. Consider
$N_G(\widetilde{v}_f)-V(MMPD)\neq\emptyset$. Let $\widetilde{v}\in
N_G(\widetilde{v}_f)-V(MMPD)$. Then, $MMPD-\{\langle
v_f,\widetilde{v}_f \rangle\}\cup\{\langle v_f,\upsilon \rangle,
\langle \widetilde{v}_f,\widetilde{v} \rangle\}$ is a
matched-paired-dominating set of $G$ having more restricted
vertices than $MMPD$, a contradiction. Thus, $\upsilon$ is not
adjacent to $v_f$ and Statement (1) holds true. Statement (2) is
clearly true. Otherwise, $MMPD\cup\{\langle \upsilon,\widetilde{v}
\rangle\}$, where $\widetilde{v}\in N_G(\upsilon)-V(MMPD)$, is a
matched-paired-dominating set of $G$ which has more restricted
vertices than $MMPD$, a contradiction.

Next, we prove Statement (3). Assume by contradiction that
$N_G(\widetilde{v}_x)-\{\upsilon\} \not\subseteq V(MMPD)$. Let
$\widetilde{v}\in (N_G(\widetilde{v}_x)-\{\upsilon\})-V(MMPD)$. By
Statement (2), $N_G(\upsilon)\subseteq V(MMPD)$. Thus,
$\widetilde{v}\not\in N_G(\upsilon)$. Then, $MMPD-\{\langle
v_x,\widetilde{v}_x \rangle\}\cup\{\langle
\widetilde{v}_x,\widetilde{v} \rangle, \langle v_x,\upsilon
\rangle\}$ is a matched-paired-dominating set of $G$ having more
restricted vertices than $MMPD$, a contradiction. Thus,
$N_G(\widetilde{v}_x)-\{\upsilon\}\subseteq V(MMPD)$.

Finally, we prove Statement (4). Assume by contradiction that
$\upsilon$ is adjacent to $v_r$. By Statement (3),
$N_G(v_f)-\{\upsilon\}\subseteq V(MMPD)$. Then, $v_f$ is dominated
by one vertex of $N_G(v_f)\cap V(MMPD)$, e.g. $v_r$. Thus,
$MMPD-\{\langle v_f,v_r \rangle\}\cup\{\langle v_r,\upsilon
\rangle\}$ is a matched-paired-dominating set of $G$ having more
restricted vertices than $MMPD$, a contradiction. Thus, $\upsilon$
is not adjacent to $v_r$.
\end{proof}

By Theorem \ref{decomposition-tree-theorem}, a decomposition tree
$DT(G)$ of a cograph $G=(V,E)$ can be constructed in
$O(|V|+|E|)$-linear time. A cograph is not connected if the root
of its corresponding decomposition tree is a $\oplus$-labeled
node. Hence, we assume that the root of the corresponding
decomposition tree is a $\otimes$-labeled node. For a node $v$ in
$DT(G)$, denote by $DT_v(G)$ the subtree of $DT(G)$ rooted at $v$,
and denote by $G_v$ the subgraph of $G$ induced by the leaves of
$DT_v(G)$. Our algorithm is sketched as follows: The algorithm is
given a decomposition tree $DT(G)$ of a cograph $G$ and a
restricted vertex set $\mathcal{R}$ in $G$. It visits nodes of
$DT(G)$ in a postorder sequence (i.e., bottom-up manner). Thus,
while visiting a node, both its children were visited. Suppose
that it is about to process internal node $v$ with $v_l$ and $v_r$
being the left and right children of $v$ in $DT(G)$, respectively.
Let $\mathcal{R}_L$ and $\mathcal{R}_R$ be the restricted vertex
sets of $G_{v_l}$ and $G_{v_r}$, respectively, such that
$|\mathcal{R}_L|\geqslant |\mathcal{R}_R|$. Let $CMPD_L$ be a
canonical matched-paired-dominating set of $G_{v_l}$ and let $V_R$
be the vertex set of $G_{v_r}$. Then, it uses $CMPD_L$ and $V_R$
to construct a canonical matched-paired-dominating set $CMPD$ of
$G_v$. If $v$ is the root of $DT(G)$, then $CMPD$ is a canonical
matched-paired-dominating set of $G$ and the algorithm terminates.
For example, let $G$ be cograph with restricted vertex set
$\mathcal{R}=\{a,f,i\}$ shown in Fig.
\ref{Fig_decomposition-tree}. Our algorithm traverses the
decomposition tree $DT(G)$ in a bottom-up manner. Suppose that it
is about to process the root $v$ with $v_l$ and $v_r$ being the
left and right children of $v$ in $DT(G)$, respectively. Then, a
canonical $(1,0,0)$-matched-paired-dominating set
$CMPD_L=\{\langle a,f \rangle\}$ of $G_{v_l}$ and the vertex set
$V_R=\{h,i\}$ of $G_{v_r}$ have been computed. The algorithm then
constructs from $CMPD_L$ and $V_R$ a canonical
$(1,1,0)$-matched-paired-dominating set $\{\langle a,f \rangle,
\langle b,i \rangle\}$ of $G_v$. In the following, we will show
how to construct such a canonical matched-paired-dominating set.

In the rest of the paper, we assume that $G=(V,E)$ is a cograph
with restricted vertex set $\mathcal{R}$ and is formed from $G_L$
and $G_R$ by either a union operation or a joint operation. We use
$V_L$ and $V_R$ to denote the vertex sets of $G_L$ and $G_R$,
respectively. In other words, $V=V_L\cup V_R$ and $V_L\cap
V_R=\emptyset$. Notice that every vertex in $V_L$ is adjacent to
all vertices in $V_R$ if $G=G_L\otimes G_R$. On the other hand, we
use $\mathcal{R}_L$ and $\mathcal{R}_R$ to denote the restricted
vertex sets of $G_L$ and $G_R$, respectively, i.e., $\mathcal{R}_L
= \mathcal{R}\cap V_L$ and $\mathcal{R}_R = \mathcal{R}\cap V_R$.

By the definition of cographs, $G_L$ or $G_R$ may contain isolated
vertices. For a graph $H$, we use $I(H)$ to denote the set of
isolated vertices in $H$. We denote by $H-I(H)$ deleting $I(H)$
from $H$. Then, $I(G_L)$ and $I(G_R)$ are the sets of isolated
vertices in $G_L$ and $G_R$, respectively. By Corollary
\ref{matched-paired-domination_exists}, $G_L-I(G_L)$ and
$G_R-I(G_R)$ have matched-paired-dominating sets if they are not
empty, and, hence, they have canonical matched-paired-dominating
sets. Then, the following lemma can be easily verified from the
definition of union operation.

\begin{lem}\label{Union-Lemma}
Assume $G=G_L\oplus G_R$ is a cograph with restricted vertex set
$\mathcal{R}$. Let $CMPD_L$ and $CMPD_R$ be canonical
matched-paired-dominating sets of $G_L-I(G_L)$ and $G_R-I(G_R)$
w.r.t. $\mathcal{R}_L-I(G_L)$ and $\mathcal{R}_R-I(G_R)$,
respectively. Then, $I(G)=I(G_L)\cup I(G_R)$ and $CMPD_L\cup
CMPD_R$ is a canonical matched-paired-dominating set of $G-I(G)$
w.r.t. $\mathcal{R}-I(G)$.
\end{lem}

From now on, we consider that $G$ is formed from $G_L$ and $G_R$
by a joint operation. First, we consider that
$\mathcal{R}=\emptyset$. Let $v_L\in V_L$ and let $v_R\in V_R$.
Obviously, $CMPD=\{\langle v_L,v_R \rangle\}$ is a
matched-paired-dominating set of $G$, and, hence, the the maximum
matched-paired-domination problem on $G$ is trivially solvable. In
the following, we assume $\mathcal{R}\neq \emptyset$. For the case
of $|\mathcal{R}_L|=|\mathcal{R}_R|$, we give the following lemma
to find a canonical matched-paired-dominating set of $G$.

\begin{lem}\label{T_L=T_R}
Assume $G=G_L\otimes G_R$ is a cograph with restricted vertex set
$\mathcal{R}$, $\mathcal{R}_L=\mathcal{R}\cap V_L$, and
$\mathcal{R}_R=\mathcal{R}\cap V_R$. If
$|\mathcal{R}_L|=|\mathcal{R}_R|$ and $|\mathcal{R}_L|>0$, then
there exists a canonical matched-paired-dominating set $CMPD$ of
$G$ w.r.t. $\mathcal{R}$ such that $V(CMPD)=\mathcal{R}$,
$|CMPD|=\frac{|\mathcal{R}|}{2}$, and $CMPD$ contains no
free-paired-edge.
\end{lem}
\begin{proof}
Let $\mathcal{R}_L=\{u_1,u_2,\cdots,u_k\}$ and
$\mathcal{R}_R=\{v_1,v_2,\cdots,v_k\}$, where
$k=\frac{|\mathcal{R}|}{2}$. By pairing $u_i$ with $v_i$ for
$1\leqslant i\leqslant k$, we obtain a
$(k,0,0)$-matched-paired-dominating set $CMPD=\{\langle u_1,v_1
\rangle, \langle u_2,v_2\rangle, \cdots, \langle u_k,v_k
\rangle\}$ of $G$ with cardinality $\frac{|\mathcal{R}|}{2}$. By
Lemma \ref{trivial_maximum_matched-paired-domination}, $CMPD$ is a
canonical matched-paired-dominating set of $G$ w.r.t.
$\mathcal{R}$ without free-paired-edges.
\end{proof}

\begin{figure}[t]
\begin{center}
\includegraphics[scale=0.7]{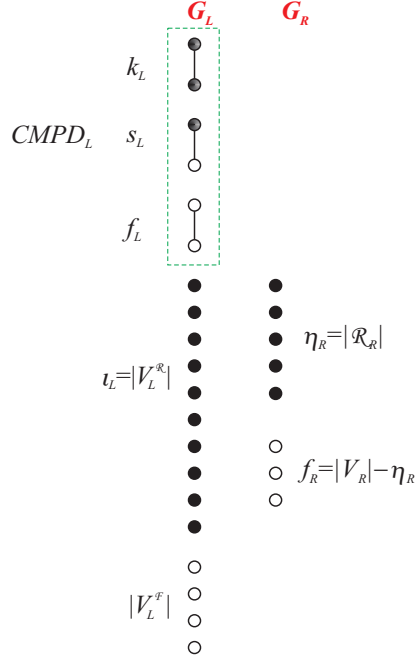}
\caption{The partition of $V_L$ and $V_R$.}\label{Fig_Partitions}
\end{center}
\end{figure}

From now on, we assume that $|\mathcal{R}_L|\neq |\mathcal{R}_R|$.
Without loss of generality, assume
$|\mathcal{R}_L|>|\mathcal{R}_R|$. Let $CMPD_L$ be a canonical
$(k_L,s_L,f_L)$-matched-paired-dominating set of $G_L-I(G_L)$
w.r.t. $\mathcal{R}_L-I(G_L)$. We first partition $V_L-V(CMPD_L)$
into two subsets $V_L^\mathcal{R}$ and $V_L^\mathcal{F}$ such that
$V_L^\mathcal{R} = \mathcal{R}_L-V(CMPD_L)$ and
$V_L^\mathcal{F}\cap \mathcal{R}_L=\emptyset$. Note that
$V_L^\mathcal{R}$ or $V_L^\mathcal{F}$ may contain isolated
vertices of $G_L$. By Statement (2) of Lemma \ref{Properties},
$N_{G_L}(v)\subseteq V(CMPD_L)$ for $v\in V_L^\mathcal{R}-I(G_L)$.
We next partition $V_R$ into two subsets $\mathcal{R}_R$ and
$V_R-\mathcal{R}_R$. The partition of $V_L$ and $V_R$ is shown in
Fig. \ref{Fig_Partitions}. For simplicity, let
$\imath_L=|V_L^\mathcal{R}|$, $\eta_R=|\mathcal{R}_R|$, and let
$f_R=|V_R|-\eta_R$. By definition,
$|\mathcal{R}_L|=2k_L+s_L+\imath_L$ and $|\mathcal{R}_R|=\eta_R$.
By assumption, $|\mathcal{R}_L|
> |\mathcal{R}_R|$. Thus, we get that
\begin{equation}\label{Assumption}
   2k_L+s_L+\imath_L > \eta_R.
\end{equation}

Considering the relation between $\imath_L$ and $\eta_R+f_R$, we
have that $\imath_L\geqslant \eta_R+f_R$ or $\imath_L <
\eta_R+f_R$. We construct from $CMPD_L$ and $V_R$ a
matched-paired-dominating set $CMPD$ of $G$ having at most one
free-paired-edge as follows:

\noindent\textbf{\textsl{Case 1:}} $\imath_L\geqslant \eta_R+f_R$.
Let $V'_L=\{u_1, u_2, \cdots, u_{|V'_L|}\}$ be a subset of
$V_L^\mathcal{R}$ such that $|V'_L|=\eta_R+f_R=|V_R|$, and let
$V_R=\{v_1, v_2, \cdots, v_{|V_R|}\}$. By pairing $u_i$ with $v_i$
for $1\leqslant i\leqslant \eta_R+f_R$, we construct a
$(k_L+\eta_R, s_L+f_R, 0)$-matched-paired-dominating set
$CMPD=K_{G_L}(CMPD_L)\cup S_{G_L}(CMPD_L)\cup_{1\leqslant
i\leqslant \eta_R+f_R} \{\langle u_i,v_i \rangle\}$.

\begin{figure}[thp]
\begin{center}
\includegraphics[scale=0.7]{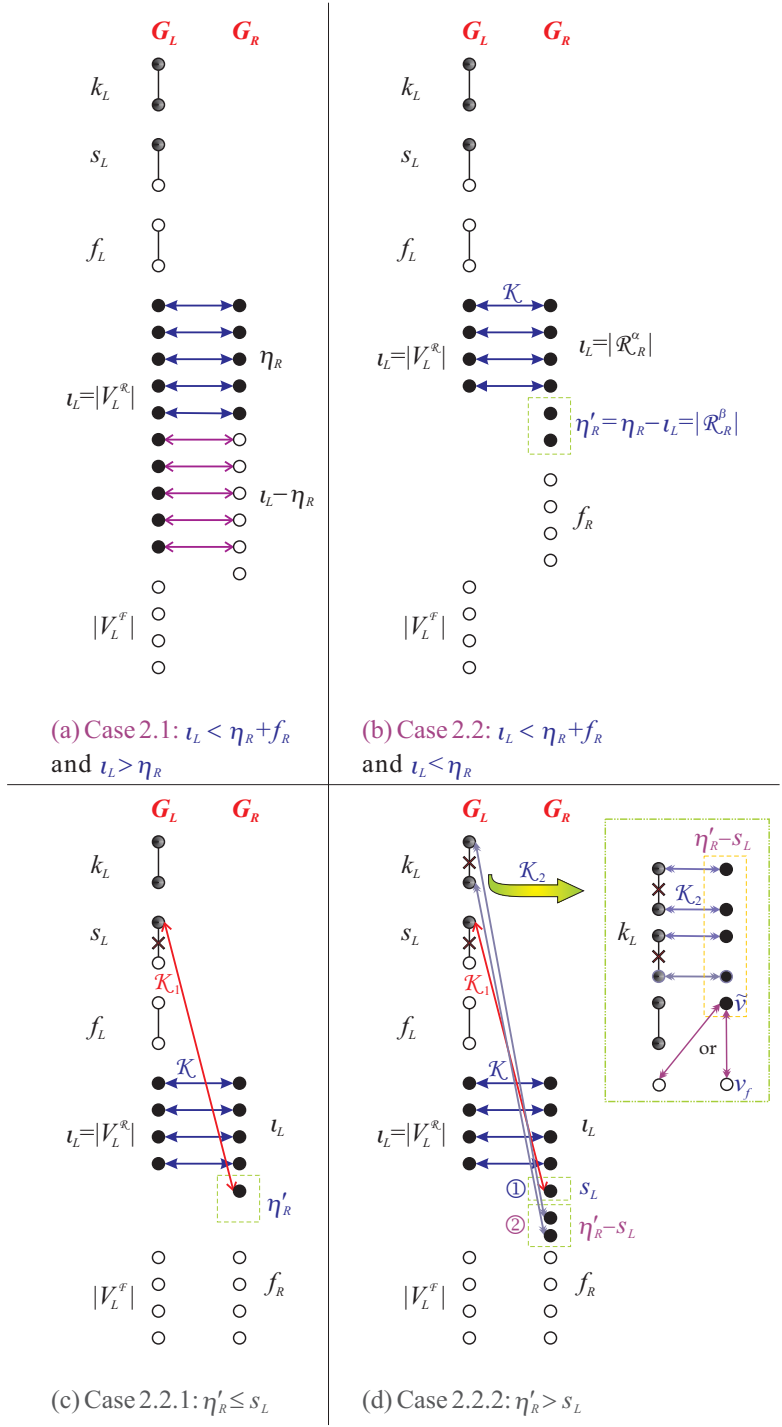}
\caption{The construction of a matched-paired-dominating set
$CMPD$ of $G$ for (a) Case 2.1, and (b)--(d) Case 2.2, where
restricted vertices are drawn by filled circles, symbol `$\times$'
denotes the destruction to one paired-edge in $CMPD_L$, and arrow
lines represent the new paired-edges for the
construction.}\label{Fig_Case2.1-2.2}
\end{center}
\end{figure}

\noindent\textbf{\textsl{Case 2:}} $\imath_L < \eta_R+f_R$. There
are three subcases:

\textbf{\textsl{Case 2.1:}} $\imath_L > \eta_R$. In this subcase,
$\eta_R < \imath_L < \eta_R+f_R$. Thus, $0 < \imath_L - \eta_R <
f_R$. Let $V_L^\mathcal{R} = \{u_1, u_2, \cdots, u_{\imath_L}\}$,
$\mathcal{R}_R=\{v_1, v_2, \cdots, v_{\eta_R}\}$, and let
$V'_R=\{v_{\eta_R+1}, v_{\eta_R+2}, \cdots, v_{\imath_L}\}$ be a
subset of $V_R-\mathcal{R}_R$ with $|V'_R|=\imath_L-\eta_R$. By
pairing $u_i$ with $v_i$ for $1\leqslant i\leqslant \imath_L$, we
obtain a $(k_L+\eta_R, s_L+\imath_L-\eta_R,
0)$-matched-paired-dominating set $CMPD=K_{G_L}(CMPD_L)\cup
S_{G_L}(CMPD_L)\cup_{1\leqslant i\leqslant \imath_L} \{\langle
u_i,v_i \rangle\}$. Then, $V(CMPD)\cap \mathcal{R} = \mathcal{R}$
and $CMPD$ contains no free-paired-edge. Thus, $CMPD$ is a
canonical matched-paired-dominating set of $G$. Fig.
\ref{Fig_Case2.1-2.2}(a) depicts the construction of $CMPD$ in the
subcase.

\textbf{\textsl{Case 2.2:}} $\imath_L < \eta_R$. By Eq.
(\ref{Assumption}), $2k_L+s_L+\imath_L > \eta_R$. Let
$\eta'_R=\eta_R-\imath_L$. Then, $2k_L+s_L > \eta_R-\imath_L =
\eta'_R > 0$. We partition $\mathcal{R}_R$ into two subsets
$\mathcal{R}_R^\alpha$ and $\mathcal{R}_R^\beta$ such that
$|\mathcal{R}_R^\alpha|=\imath_L$ and
$|\mathcal{R}_R^\beta|=\eta'_R=\eta_R-\imath_L$. By pairing every
vertex in $\mathcal{R}_R^\alpha$ with a vertex in
$V_L^\mathcal{R}$, we obtain a set $\mathcal{K}$ of $\imath_L$
full-paired-edges shown in Fig. \ref{Fig_Case2.1-2.2}(b). We then
consider the following two subcases:

\hspace{1cm}\textbf{\textsl{Case 2.2.1:}}
$\eta_R-\imath_L=\eta'_R\leqslant s_L$. We first partition
$S_{G_L}(CMPD_L)$ into two subsets $\mathcal{S}_1$ and
$\mathcal{S}_2$ such that $|\mathcal{S}_1|=\eta'_R$ and
$|\mathcal{S}_2|=s_L-\eta'_R$. Let $u_1, u_2, \cdots, u_{\eta'_R}$
be the restricted vertices in $V(\mathcal{S}_1)$ and let
$\mathcal{R}_R^\beta=\{v_1, v_2, \cdots, v_{\eta'_R}\}$. By
pairing $u_i$ with $v_i$ for $1\leqslant i\leqslant \eta'_R$, we
obtain a set $\mathcal{K}_1$ of $\eta'_R$ full-paired-edges. Let
$CMPD=K_{G_L}(CMPD_L)\cup \mathcal{K}\cup \mathcal{K}_1 \cup
\mathcal{S}_2$. Then, $CMPD$ is a $(k_L+\eta_R, s_L-\eta'_R,
0)$-matched-paired-dominating set of $G$. Since $V(CMPD)\cap
\mathcal{R} = \mathcal{R}$, $CMPD$ is a maximum
matched-paired-dominating set of $G$. Thus, $CMPD$ is a canonical
matched-paired-dominating set of $G$. The construction of $CMPD$
is shown in Fig. \ref{Fig_Case2.1-2.2}(c).

\hspace{1cm}\textbf{\textsl{Case 2.2.2:}} $\eta_R-\imath_L=\eta'_R
> s_L$. We first partition $\mathcal{R}_R^\beta$ into two subsets
$\mathcal{R}_R^\textrm{a}$ and $\mathcal{R}_R^\textrm{b}$ such
that $|\mathcal{R}_R^\textrm{a}|=s_L$ and
$|\mathcal{R}_R^\textrm{b}|=\eta'_R-s_L$. Let $u_1, u_2, \cdots,
u_{s_L}$ be the restricted vertices in $V(S_{G_L}(CMPD_L))$ and
let $\mathcal{R}_R^\textrm{a}=\{v_1, v_2, \cdots, v_{s_L}\}$. By
pairing $u_i$ with $v_i$ for $1\leqslant i\leqslant s_L$, we
obtain a set $\mathcal{K}_1$ of $s_L$ full-paired-edges. Suppose
that $|\mathcal{R}_R^\textrm{b}|=\eta'_R-s_L$ is even. We first
partition $K_{G_L}(CMPD_L)$ into two subsets $\mathcal{K}_{L_1}$
and $\mathcal{K}_{L_2}$ such that $\mathcal{K}_{L_1}$ contains
$\frac{|\mathcal{R}_R^\textrm{b}|}{2}$ full-paired-edges, i.e.,
$V(\mathcal{K}_{L_1})$ contains $|\mathcal{R}_R^\textrm{b}|$
restricted vertices. By pairing every vertex of
$V(\mathcal{K}_{L_1})$ with a vertex in
$\mathcal{R}_R^\textrm{b}$, we obtain a set $\mathcal{K}_2$ of
$|\mathcal{R}_R^\textrm{b}|$ full-paired-edges. Let
$CMPD=\mathcal{K}\cup\mathcal{K}_1\cup\mathcal{K}_2
\cup\mathcal{K}_{L_2}$. Then, $CMPD$ is a
$(k_L+\frac{\imath_L+s_L+\eta_R}{2}, 0,
0)$-matched-paired-dominating set of $G$. We can see that $CMPD$
is a $(\frac{|\mathcal{R}|}{2}, 0, 0)$-matched-paired-dominating
set of $G$, $V(CMPD)\cap \mathcal{R} = \mathcal{R}$, and that
$CMPD$ contains no free-paired-edge. On the other hand, suppose
that $|\mathcal{R}_R^\textrm{b}|=\eta'_R-s_L$ is odd. Then,
$|\mathcal{R}|$ is odd. We pick from $G$ a restricted vertex
$\widetilde{v}$ and a free vertex $v_f$ to form a semi-paired-edge
as follows: Case i, $V_L-\mathcal{R}_L\neq\emptyset$. Let
$\widetilde{v}\in \mathcal{R}_R^\textrm{b}$ and let $v_f\in
V_L-\mathcal{R}_L$. Case ii, $V_L-\mathcal{R}_L=\emptyset$ and
$f_R>|I(G_R)-\mathcal{R}_R|$. Let $\widetilde{v}$ be a restricted
vertex in $\mathcal{R}_R^\textrm{b}$ such that $\widetilde{v}$ is
adjacent to one free vertex $v_f$ in $V_R-(\mathcal{R}_R\cup
I(G_R))$. Case iii, $V_L-\mathcal{R}_L=\emptyset$ and
$f_R=|I(G_R)-\mathcal{R}_R|\neq 0$. Let $\langle \widetilde{v},
v_L\rangle$ be a full-paired-edge in $K_{G_L}(CMPD_L)$ and let
$v_f$ be a free vertex in $I(G_R)-\mathcal{R}_R$. For case of
$V_L-\mathcal{R}_L=\emptyset$ and $f_R=|I(G_R)-\mathcal{R}_R|=0$,
we have that $|V|=|\mathcal{R}|$ is odd and a
$(\lfloor\frac{|\mathcal{R}|}{2}\rfloor, 0,
0)$-matched-paired-dominating set $CMPD$ of $G$ can be easily
constructed from $K_{G_L}(CMPD_L)$ and $\mathcal{R}_R$. By
Statement (2) of Lemma
\ref{trivial_maximum_matched-paired-domination}, $CMPD$ is a
canonical matched-paired-dominating set of $G$. Now, suppose
$\widetilde{v}$ and $v_f$ exist. Let $\mathcal{S}=\{\langle
\widetilde{v},v_f \rangle\}$. If $\widetilde{v}\not\in
\mathcal{R}_L$, then let $\widetilde{\mathcal{K}}=\emptyset$ and
$\mathcal{R}_R^\textrm{b}=\mathcal{R}_R^\textrm{b}-\{\widetilde{v}\}$;
otherwise, let $\langle \widetilde{v},v_L \rangle\in
K_{G_L}(CMPD_L)$, $K_{G_L}(CMPD_L)=K_{G_L}(CMPD_L)-\{\langle
\widetilde{v},v_L \rangle\}$, $v_R\in \mathcal{R}_R^\textrm{b}$,
$\mathcal{R}_R^\textrm{b}=\mathcal{R}_R^\textrm{b}-\{v_R\}$, and
let $\widetilde{\mathcal{K}}=\{\langle v_L,v_R \rangle\}$. Then,
$|\mathcal{R}_R^\textrm{b}|$ becomes even. We then partition
$K_{G_L}(CMPD_L)$ into two subsets $\mathcal{K}_{L_1}$ and
$\mathcal{K}_{L_2}$ such that $\mathcal{K}_{L_1}$ contains
$\frac{|\mathcal{R}_R^\textrm{b}|}{2}$ full-paired-edges. By
pairing every vertex of $V(\mathcal{K}_{L_1})$ with a vertex in
$\mathcal{R}_R^\textrm{b}$, we obtain a set $\mathcal{K}_2$ of
$|\mathcal{R}_R^\textrm{b}|$ full-paired-edges. Let
$CMPD=\mathcal{K}\cup\mathcal{K}_1\cup\mathcal{K}_2
\cup\mathcal{K}_{L_2}\cup\widetilde{\mathcal{K}}\cup\mathcal{S}$.
Then, $CMPD$ is a
$(k_L+\lfloor\frac{\imath_L+s_L+\eta_R}{2}\rfloor, 1,
0)$-matched-paired-dominating set of $G$. We can see that $CMPD$
is a $(\lfloor\frac{|\mathcal{R}|}{2}\rfloor, 1,
0)$-matched-paired-dominating set of $G$. By Statement (1) of
Lemma \ref{trivial_maximum_matched-paired-domination}, $CMPD$ is a
canonical matched-paired-dominating set of $G$. The construction
of $CMPD$ is shown in Fig. \ref{Fig_Case2.1-2.2}(d).

\begin{figure}[thp]
\begin{center}
\includegraphics[scale=0.7]{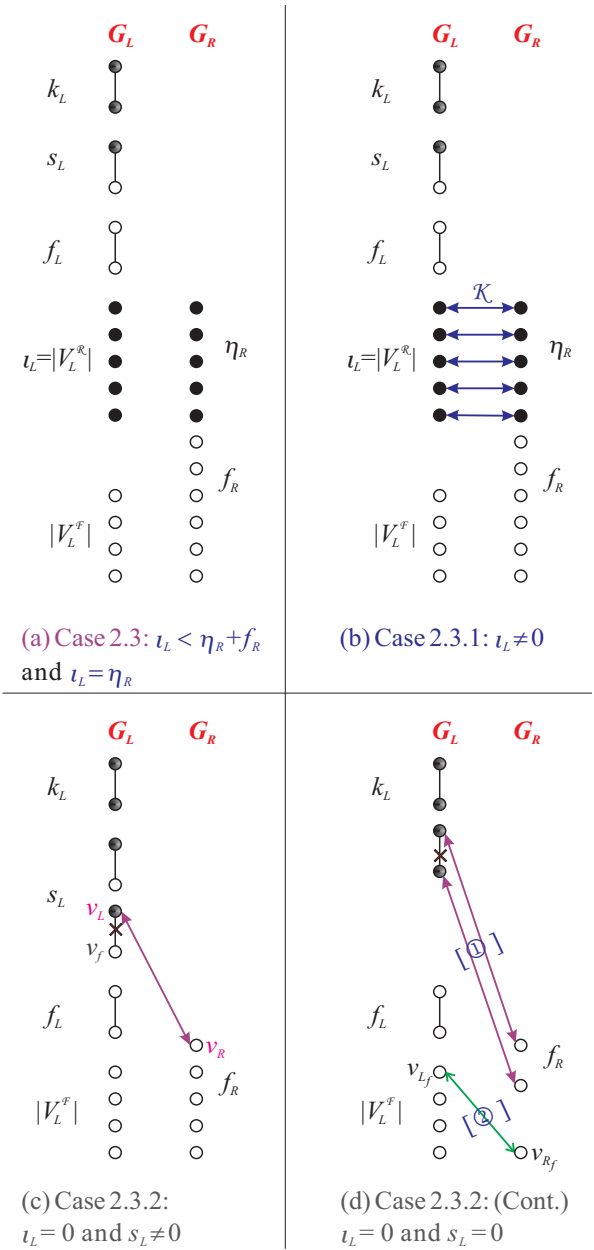}
\caption{The construction of a matched-paired-dominating set
$CMPD$ of $G$ for Case 2.3, where (a) the partition of $V_L$ and
$V_R$ for the case, (b) the construction of a
matched-paired-dominating set for case of $\imath_L\neq 0$, and
(c)--(d) the construction of a matched-paired-dominating set for
case of $\imath_L=0$. Note that restricted vertices are drawn by
filled circles, symbol `$\times$' denotes the destruction to one
paired-edge in $CMPD_L$, and arrow lines represent the new
paired-edges for the construction.}\label{Fig_Case2.3}
\end{center}
\end{figure}

\textbf{\textsl{Case 2.3:}} $\imath_L = \eta_R$. In this subcase,
$\imath_L = \eta_R < \eta_R+f_R$. Consider the following two
subcases:

\hspace{1cm}\textbf{\textsl{Case 2.3.1:}} $\imath_L \neq 0$. Let
$V_L^\mathcal{R}=\{u_1, u_2, \cdots, u_{\imath_L}\}$ be the
restricted vertex set of $V_L-V(CMPD_L)$ and let
$\mathcal{R}_R=\{v_1, v_2, \cdots, v_{\imath_L}\}$. By paring
$u_i$ with $v_i$ for $1\leqslant i\leqslant \imath_L$, we get a
set $\mathcal{K}$ of $\imath_L$ full-paired-edges. Let
$CMPD=K_{G_L}(CMPD_L)\cup S_{G_L}(CMPD_L)\cup \mathcal{K}$. Then,
$CMPD$ is a $(k_L+\imath_L, s_L, 0)$-matched-paired-dominating set
of $G$. We can see that $CMPD$ is a maximum
matched-paired-dominating set of $G$ without free-paired-edges.
Thus, $CMPD$ is a canonical matched-paired-dominating set of $G$.
Fig. \ref{Fig_Case2.3}(b) shows the construction of $CMPD$ in this
subcase.

\hspace{1cm}\textbf{\textsl{Case 2.3.2:}} $\imath_L = 0$. First,
we consider that $s_L\neq 0$. Let $\langle v_L, v_f \rangle$ be a
semi-paired-edge, with restricted vertex $v_L$, in
$S_{G_L}(CMPD_L)$, and let $v_R$ be a free vertex in $V_R$. Let
$\mathcal{S}=\{\langle v_L, v_R \rangle\}$. Then,
$CMPD=K_{G_L}(CMPD_L)\cup S_{G_L}(CMPD_L)-\{\langle v_L, v_f
\rangle\}\cup\mathcal{S}$ is a $(k_L, s_L,
0)$-matched-paired-dominating set of $G$ such that $V(CMPD)\cap
\mathcal{R} = \mathcal{R}$. It is easy to see that $CMPD$ is a
canonical matched-paired-dominating set of $G$. Fig.
\ref{Fig_Case2.3}(c) shows the construction of $CMPD$ in case of
$\imath_L = \eta_R = 0$ and $s_L \neq 0$. On the other hand, we
consider that $s_L=0$. If $V(K_{G_L}(CMPD_L))$ is a dominating set
of $G_L$, i.e., $f_L=0$ and $I(G_L)=\emptyset$, then
$CMPD=K_{G_L}(CMPD_L)$ is clearly a canonical
matched-paired-dominating set of $G$. Suppose that $f_L\neq 0$ or
$I(G_L)\neq\emptyset$. Consider that $f_R\geqslant 2$. Let
$v_{f_1}$ and $v_{f_2}$ be two free vertices in $V_R$, and let
$\langle v_{r_1},v_{r_2} \rangle$ be a full-paired-edge in
$K_{G_L}(CMPD_L)$. Then, $CMPD = K_{G_L}(CMPD_L)-\{\langle
v_{r_1},v_{r_2} \rangle\}\cup \{\langle v_{f_1},v_{r_1} \rangle,
\langle v_{f_2},v_{r_2} \rangle\}$ is a canonical $(k_L-1, 2,
0)$-matched-paired-dominating set of $G$ with that $V(CMPD)\cap
\mathcal{R} = \mathcal{R}$. Next, consider that $f_R=1$. Let
$v_{R_f}$ be the only vertex in $V_R$. Consider the following
cases: Case i, there exists one restricted vertex
$\widetilde{v}_L$ in $\mathcal{R}_L$ such that
$N_{G_L}(\widetilde{v}_L)\not\subseteq \mathcal{R}_L$. Let
$v_{L_f}\in N_{G_L}(\widetilde{v}_L)-\mathcal{R}_L$ and let
$\langle \widetilde{v}_L,v_L\rangle$ be a full-paired-edge in
$K_{G_L}(CMPD_L)$. Let $CMPD=K_{G_L}(CMPD_L)-\{\langle
\widetilde{v}_L,v_L \rangle\}\cup\{\langle v_L,v_{R_f} \rangle,
\langle \widetilde{v}_L,v_{L_f} \rangle\}$. Then, $CMPD$ is a
canonical $(k_L-1, 2, 0)$-matched-paired-dominating set of $G$
with that $V(CMPD)\cap \mathcal{R} = \mathcal{R}$. Case ii,
$N_{G_L}(\widetilde{v}_L)\subseteq \mathcal{R}_L$ for each
$\widetilde{v}_L\in\mathcal{R}_L$. Let $v_{L_f}\in
V_L-\mathcal{R}_L$ and let $CMPD=K_{G_L}(CMPD_L)\cup\{\langle
v_{L_f},v_{R_f} \rangle\}$. Then, $CMPD$ is a $(k_L, 0,
1)$-matched-paired-dominating set of $G$. We can see that if
$N_{G_L}(\widetilde{v}_L)\subseteq \mathcal{R}_L$ for each
$\widetilde{v}_L\in\mathcal{R}_L$, then a free-paired-edge is
necessary for constructing a maximum matched-paired-dominating set
of $G$. Thus, $CMPD$ is a canonical matched-paired-dominating set
of $G$. Fig. \ref{Fig_Case2.3}(d) depicts the construction of
$CMPD$ in case of $\imath_L = \eta_R = 0$ and $s_L =
0$.\hfill{$\square$}

It follows from the above constructions and arguments that our
constructed matched-paired-dominating set $CMPD$ for case of
$\imath_L < \eta_R+f_R$ (Case 2) is a canonical
matched-paired-dominating set of $G$. The remnant is to prove that
the constructed matched-paired-dominating set $CMPD$ for case of
$\imath_L \geqslant \eta_R+f_R$ (Case 1) is a canonical
matched-paired-dominating set of $G$. The following lemma shows
the result.

\begin{lem}\label{Case_1}
Assume $G=G_L\otimes G_R$ is a cograph with restricted vertex set
$\mathcal{R}$, $\mathcal{R}_L=\mathcal{R}\cap V_L$,
$\mathcal{R}_R=\mathcal{R}\cap V_R$, and
$|\mathcal{R}_L|>|\mathcal{R}_R|$. Let $CMPD_L$ be a canonical
$(k_L,s_L,f_L)$-matched-paired-dominating set of $G_L-I(G_L)$,
$\imath_L=|\mathcal{R}_L-V(CMPD_L)|$, $\eta_R=|\mathcal{R}_R|$,
and let $f_R=|V_R|-\eta_R$. If $\imath_L\geqslant \eta_R+f_R$,
then the constructed $(k_L+\eta_R, s_L+f_R,
0)$-matched-paired-dominating set $CMPD$ is a canonical
matched-paired-dominating set of $G$.
\end{lem}
\begin{proof}
In case of $\imath_L\geqslant \eta_R+f_R$, the construction of
$CMPD$ is shown in Fig. \ref{Fig_Case1}(a). A paired-edge in a
matched-paired-dominating set of $G$ is called \textit{mixed} if
one of its vertices is in $V_L$ and the other is in $V_R$. Suppose
that $MMPD$ is a maximum matched-paired-dominating set of $G$ with
the least free-paired-edges. That is, $MMPD$ is a canonical
matched-paired-dominating set of $G$. We may assume that $MMPD$ is
chosen such that the number of mixed paired-edges is maximal.
Denote by $MMPD_{|G_L}$ (resp. $MMPD_{|G_R}$) a restriction of
$MMPD$ to $G_L$ (resp. $G_R$). The set of mixed paired-edges of
$MMPD$ is partitioned into four subsets $K,S_1,S_2,F$ such that
$K$ contains all mixed full-paired-edges, $S_1$ contains all mixed
semi-paired-edges with restricted vertices being in $V_L$, $S_2$
contains all mixed semi-paired-edges with restricted vertices
being in $V_R$, and $F$ contains all mixed free-paired-edges. The
set of paired-edges of $MMPD_{|G_L}$ (resp. $MMPD_{|G_R}$) is
partitioned into three subsets $K_L$, $S_L$, $F_L$ (resp. $K_R$,
$S_R$, and $F_R$) containing full-paired-edges, semi-paired-edges,
and free-paired-edges, respectively. Let
$I_L=\mathcal{R}_L-V(MMPD)$ and $I_R=\mathcal{R}_R-V(MMPD)$. For
simplicity, let $|K|=k$, $|S_1|=s_1$, $|S_2|=s_2$, $|F|=f$,
$|K_L|=k'_L$, $|S_L|=s'_L$, $|F_L|=f'_L$, $|K_R|=k'_R$,
$|S_R|=s'_R$, $|F_R|=f'_R$, $|I_L|=i'_L$, and let $|I_R|=i'_R$.
The possible paired-edges in $MMPD$ are shown in Fig.
\ref{Fig_Case1}(b). Since $|\mathcal{R}_L| > |\mathcal{R}_R|$ and
$MMPD$ is a canonical matched-paired-dominating set of $G$,
$f\leqslant 1$ and at least one of $i'_L$ and $i'_R$ equals to 0.

\begin{figure}[t]
\begin{center}
\includegraphics[scale=0.7]{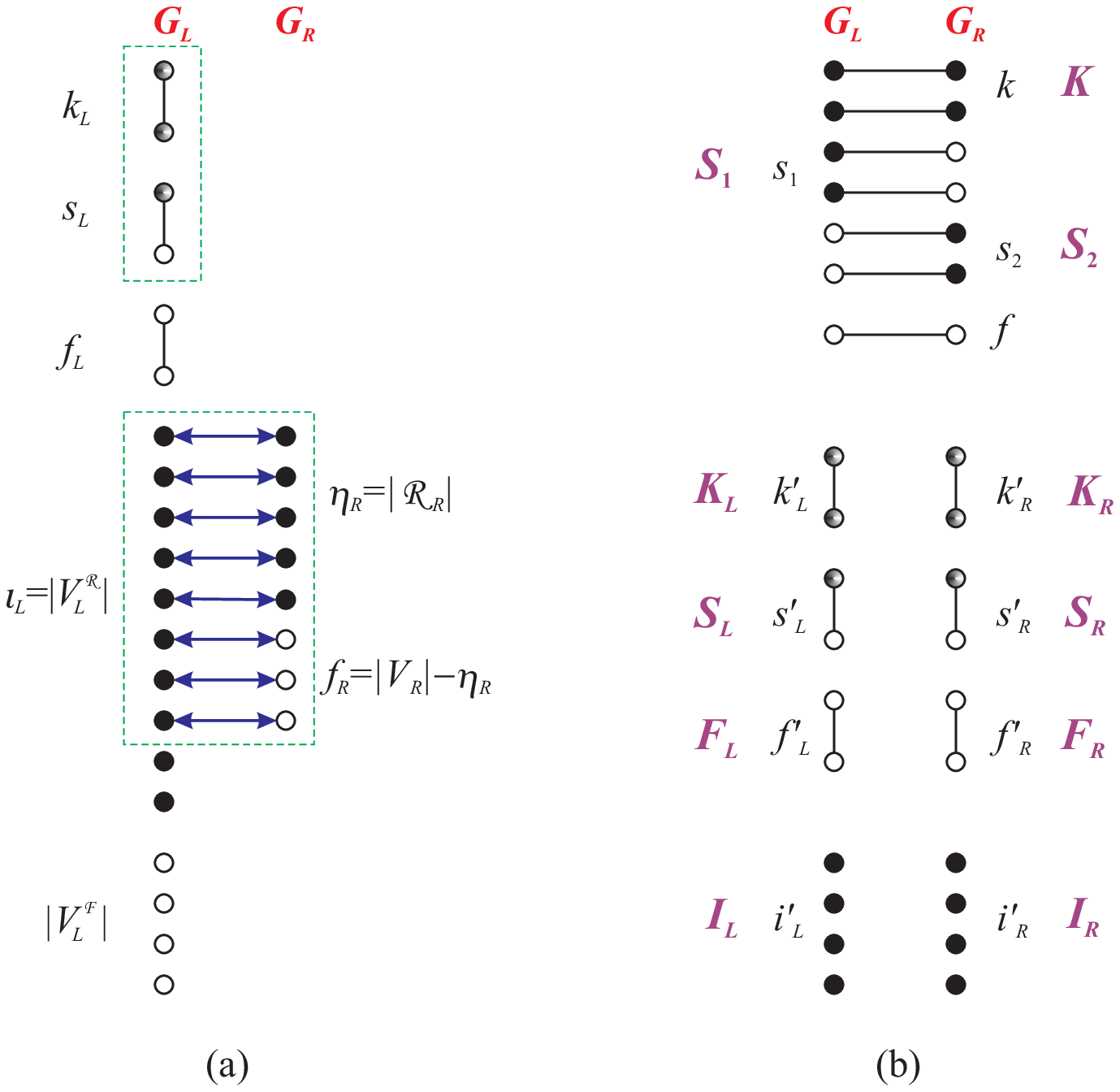}
\caption{(a) The construction of a matched-paired-dominating set
$CMPD$ of $G$ for case of $\imath_L\geqslant \eta_R+f_R$, and (b)
the possible paired-edges in a canonical matched-paired-dominating
set $MMPD$ of $G$ with the largest number of mixed paired-edges,
where restricted vertices are drawn by filled circles and arrow
lines represent the new paired-edges for the
construction.}\label{Fig_Case1}
\end{center}
\end{figure}

We first prove Claim (1) that $2k_L+s_L\geqslant 2k'_L+s'_L$. We
prove it by constructing from $MMPD_{|G_L}$ a
matched-paired-dominating set $MMPD_L$ of $G_L-I(G_L)$ such that
$|V(MMPD_L)\cap(\mathcal{R}_L-I(G_L))| \geqslant 2k'_L + s'_L$.
The construction is as follows: Initially, let $MMPD_L=K_L\cup
S_L$. Let $V'_L = V_L-I(G_L)-V(MMPD_L)$. For $v'_L\in V'_L$, if
$v'_L$ is not dominated by $V(MMPD_L)$, then let $v''_L\in
N_{G_L}(v'_L)-V(MMPD_L)$, $MMPD_L = MMPD_L\cup \{\langle v'_L,
v''_L\rangle\}$, and let $V'_L = V'_L- \{v'_L, v''_L\}$;
otherwise, let $V'_L = V'_L- \{v'_L\}$. Since $v'_L$ is not an
isolated vertex in $G_L$, $v''_L$ does exist if $v'_L$ is not
dominated by $V(MMPD_L)$. Repeat the above process until
$V'_L=\emptyset$. Then, $MMPD_L$ is a matched-paired-dominating
set of $G_L-I(G_L)$ satisfying that
$|V(MMPD_L)\cap(\mathcal{R}_L-I(G_L))|\geqslant 2k'_L + s'_L$.
Since $CMPD_L$ is a maximum $(k_L, s_L,
f_L)$-matched-paired-dominating set of $G_L-I(G_L)$,
$|V(CMPD_L)\cap (\mathcal{R}_L-I(G_L))|= 2k_L+s_L\geqslant
|V(MMPD_L)\cap (\mathcal{R}_L-I(G_L))|\geqslant 2k'_L + s'_L$.

Next, we prove Claim (2) that $i'_R=s_2=k'_R=s'_R=0$. We first
show that $i'_R=0$. Assume by contradiction that $i'_R\neq 0$.
Then, $i'_L=0$. By Statements (1) and (4) of Lemma
\ref{Properties}, $s_1=f=s'_L=f'_L = 0$. By assumption,
$|\mathcal{R}_L|=k+2k'_L > |\mathcal{R}_R|=k+s_2+2k'_R+s'_R+i'_R$.
Thus, $2k'_L > s_2+2k'_R+s'_R+i'_R \geqslant 1$ and, hence,
$k'_L\geqslant 1$. Let $v_R$ be a restricted vertex in $I_R$ and
let $\langle v_L,v'_L \rangle$ be a full-paired-edge in $K_L$. By
pairing $v_L$ with $v_R$ and all the other paired-edges stay the
same, we obtain a maximum matched-paired-dominating set $MMPD'$ of
$G$ having more mixed paired-edges than $MMPD$, a contradiction.
Thus, $i'_R=0$. We then prove that $s_2=k'_R=s'_R=0$. Assume by
contradiction that $s_2+k'_R+s'_R\neq0$. By assumption,
$|\mathcal{R}_L| > |\mathcal{R}_R|$. Then, $k+s_1+2k'_L+s'_L+i'_L
> k+s_2+2k'_R+s'_R$. Thus, $2k'_L+s_1+s'_L+i'_L >
2k'_R+s_2+s'_R$. Let $R$ be the set of restricted vertices in
$S_2\cup K_R\cup S_R$. Then, $|R|=2k'_R+s_2+s'_R$. Suppose that
$2k'_L < 2k'_R+s_2+s'_R$. Let $L$ be a subset of restricted
vertices in $S_1\cup S_L\cup I_L$ such that
$|L|=(2k'_R+s_2+s'_R)-2k'_L$. Then,
$|L|+|V(K_L)|=2k'_R+s_2+s'_R=|R|$. By pairing every vertex in $R$
with one restricted vertex of $V(K_L)\cup L$ and all the other
paired-edges stay the same, we obtain a maximum
matched-paired-dominating set $MMPD'$ of $G$ having more mixed
paired-edges than $MMPD$, a contradiction. In the following,
suppose that $2k'_L\geqslant 2k'_R+s_2+s'_R$. Consider the
following cases:

\noindent \textit{Case 1: $2k'_R+s_2+s'_R$ is even.} Let
$K_L=K_L^a\cup K_L^b$ such that $K_L^a\cap K_L^b=\emptyset$ and
$|K_L^a|=\frac{2k'_R+s_2+s'_R}{2}$. By pairing every vertex in $R$
with one restricted vertex of $V(K_L^a)$ and all the other
paired-edges stay the same, we obtain a maximum
matched-paired-dominating set $MMPD'$ of $G$ having more mixed
paired-edges than $MMPD$, a contradiction.

\noindent \textit{Case 2: $2k'_R+s_2+s'_R$ is odd.} Let
$K_L=K_L^a\cup K_L^b$ such that $K_L^a\cap K_L^b=\emptyset$ and
$|K_L^a|=\lfloor\frac{2k'_R+s_2+s'_R}{2}\rfloor$. Consider the
following subcases:

\textit{Case 2.1: $s_1+s'_L+i'_L\neq 0$.} Let $v_L$ be a
restricted vertex in $S_1\cup S_L\cup I_L$. Let $v_R$ be a vertex
in $R$ and let $R'=R-\{v_R\}$. Then,
$\frac{|R'|}{2}=\lfloor\frac{2k'_R+s_2+s'_R}{2}\rfloor$. By
pairing $v_R$ with $v_L$, pairing every vertex in $R'$ with one
restricted vertex of $V(K_L^a)$, and all the other paired-edges
stay the same, we obtain a maximum matched-paired-dominating set
$MMPD'$ of $G$ having more mixed paired-edges than $MMPD$, a
contradiction.

\textit{Case 2.2: $s_1+s'_L+i'_L=0$.} Suppose that $s'_R\neq 0$.
Let $\langle v_R,v_f\rangle$ be a semi-paired-edge, with
restricted vertex $v_R$, in $S_R$. Let $R'=R-\{v_R\}$ and let
$\langle v_L,v'_L \rangle$ be a full-paired-edge in $K_L^b$. By
pairing $v_R$ with $v_L$, pairing $v_f$ with $v'_L$, pairing every
vertex in $R'$ with one restricted vertex of $V(K_L^a)$, and all
the other paired-edges stay the same, we obtain a maximum
matched-paired-dominating set $MMPD'$ of $G$ having more mixed
paired-edges than $MMPD$, a contradiction. On the other hand,
suppose that $s'_R=0$. Then, $s_2$ is odd. We prove $f_R\neq 0$.
Assume by contradiction that $f_R=0$. By assumption of the lemma,
$\imath_L\geqslant\eta_R=|V_R|=k+s_2+2k'_R$. Since $s_2>0$,
$\imath_L\geqslant k+1$ and, hence, $\imath_L-k>0$. Then,
$|\mathcal{R}_L|=2k_L+s_L+\imath_L=k+2k'_L$. Consequently,
$(2k_L+s_L)-2k'_L=k-\imath_L<0$. It contradicts that
$(2k_L+s_L)-2k'_L\geqslant 0$ by Claim (1). Thus, $f_R\neq 0$. Let
$v_f$ be a free vertex in $V_R$, $\langle v_R,v'_f\rangle$ be a
semi-paired-edge in $S_2$ such that $v_R\in \mathcal{R}_R$,
$R'=R-\{v_R\}$, and let $\langle v_L,v'_L \rangle$ be a
full-paired-edge in $K_L^b$. By pairing $v_R$ with $v_L$, pairing
$v_f$ with $v'_L$, pairing every vertex in $R'$ with one
restricted vertex of $V(K_L^a)$, and all the other paired-edges
stay the same, we obtain a maximum matched-paired-dominating set
$MMPD'$ of $G$ having more mixed paired-edges than $MMPD$, a
contradiction.

It follows from the above arguments that Claim (2) holds true;
i.e., $i'_R=s_2=k'_R=s'_R=0$. Thus, $k=\eta_R$. Suppose that
$i'_L\neq 0$. Then, $f=f'_R=0$ by Statement (1) of Lemma
\ref{Properties}. Assume by contradiction that
$\widehat{f}_R=f_R-s_1\neq 0$. Then,
$|\mathcal{R}_L|=k+s_1+2k'_L+s'_L+i'_L=2k_L+s_L+\imath_L\geqslant
2k_L+s_L+\eta_R+s_1+\widehat{f}_R=2k_L+s_L+k+s_1+\widehat{f}_R$.
By Claim (1), $2k_L+s_L\geqslant 2k'_L+s'_L$. Thus, $i'_L\geqslant
\widehat{f}_R$. By pairing every free vertex in $V_R-V(MMPD)$ with
one restricted vertex in $I_L$ and all the other paired-edges stay
the same, we obtain a matched-paired-dominating set of $G$ having
more restricted vertices than $MMPD$, a contradiction. Thus,
$f_R=s_1$. On the other hand, suppose that $i'_L=0$. By Claim (1),
$2k_L+s_L\geqslant 2k'_L+s'_L$. By assumption of the lemma,
$\imath_L\geqslant |V_R|=k+s_1+f_R-s_1$. Then,
$|\mathcal{R}_L|=k+s_1+2k'_L+s'_L=2k_L+s_L+\imath_L\geqslant
2k_L+s_L+k+s_1+f_R-s_1$. Thus, $2k'_L+s'_L\geqslant
2k_L+s_L+f_R-s_1$. Since $2k_L+s_L\geqslant 2k'_L+s'_L$ by Claim
(1), $f_R-s_1=0$. Thus, $f_R=s_1$. Consequently, $k=\eta_R$ and
$s_1=f_R$. We can see that $|V(CMPD)\cap \mathcal{R}| =
2k_L+s_L+2\eta_R+f_R=2k_L+s_L+2k+s_1$ and $|V(MMPD)\cap
\mathcal{R}| = 2k'_L+s'_L+2k+s_1$. Thus, $|V(CMPD)\cap
\mathcal{R}|-|V(MMPD)\cap
\mathcal{R}|=(2k_L+s_L)-(2k'_L+s'_L)\geqslant 0$ by Claim (1).
That is, the constructed matched-paired-dominating set $CMPD$ is a
maximum matched-paired-dominating set of $G$. In addition, $CMPD$
contains no free-paired-edge. Thus, the constructed $(k_L+\eta_R,
s_L+f_R, 0)$-matched-paired-dominating set $CMPD$ is a canonical
matched-paired-dominating set of $G$.
\end{proof}

It follows from Lemma \ref{Case_1} that our constructed
matched-paired-dominating set $CMPD$ is a canonical
matched-paired-dominating set of $G$ w.r.t. $\mathcal{R}$. Now, we
will analyze the time complexity for constructing $CMPD$. For case
of $\imath_L\geqslant \eta_R+f_R$ shown in Fig.
\ref{Fig_Case1}(a), $CMPD$ is constructed in $O(|V_R|)$ time,
where $|V_R|\leqslant |\mathcal{R}_L|$. Consider that
$\imath_L<\eta_R+f_R$. For case of $\eta_R<\imath_L$ shown in Fig.
\ref{Fig_Case2.1-2.2}(a), $CMPD$ is constructed in $O(\imath_L)$
time, where $\imath_L\leqslant |\mathcal{R}_L|$. For case of
$\eta_R>\imath_L$ shown in Fig. \ref{Fig_Case2.1-2.2}(b)--(d),
$CMPD$ can be easily constructed in $O(|\mathcal{R}_R|)$ time,
where $|\mathcal{R}_R|\leqslant |\mathcal{R}_L|$. On the other
hand, for case of $\imath_L=\eta_R$ shown in Fig.
\ref{Fig_Case2.3}, $CMPD$ can be constructed in
$O(|\mathcal{R}_R|)$ time, where $|\mathcal{R}_R|\leqslant
|\mathcal{R}_L|$. It follows from the above arguments that
constructing a canonical matched-paired-dominating set $CMPD$ of
$G$ runs in $O(|\mathcal{R}_L|)$ time. Let
$\widehat{E}_{LR}=\{uv|\forall u\in V_L$ and $\forall v\in V_R\}$.
Then, $|\mathcal{R}_L|\leqslant |\widehat{E}_{LR}|$. Hence, a
canonical matched-paired-dominating set $CMPD$ of $G$ can be
computed in $O(|\widehat{E}_{LR}|)$ time.

It follows from the above analysis that given a decomposition tree
of a cograph $G=(V,E)$ and a restricted vertex set
$\mathcal{R}\subseteq V$, a canonical matched-paired-dominating
set of $G$ w.r.t. $\mathcal{R}$ can be constructed in
$O(|V|+|E|)$-linear time. Thus, we conclude the following theorem.

\begin{thm}\label{MainResult}
Given a cograph $G=(V,E)$ with restricted vertex set
$\mathcal{R}$, the maximum matched-paired-domination problem can
be solved in $O(|V|+|E|)$-linear time.
\end{thm}

\section{Concluding Remarks}
The paired-domination problem can be applied to allocate guards on
vertices such that these guards protect every vertex, each guard
is assigned another adjacent one, and they are designed as backup
for each other. However, some vertices may play more important
role (for example, important facilities are placed on these
vertices) and, hence, they are placed by guards for instant
protection possible. Motivated by the issue we propose a
generalization of the paired-domination problem, namely, the
maximum matched-paired-domination problem. We then solve the
maximum matched-paired-domination problem on cographs in linear
time. A future work will be to extend our technique to solve the
maximum matched-paired-domination problem on some special classes
of graphs, such as trees, block graphs, Ptolemaic graphs and
distance-hereditary graphs.

\newpage
\section*{\Large List of Symbols}

\fontsize{11}{13pt}\selectfont
\begin{itemize}
  \item[1.] $N_G(v)$, $N_G[v]$: $N_G(v)$ is the \textit{open neighborhood} of a vertex $v$ in a graph
  $G=(V,E)$ and is defined to be $\{u\in V |uv\in E\}$. $N_G[v]$ is the \textit{closed neighborhood} of
  $v$ and is defined to be $N_G(v)\cup\{v\}$.
  \item[2.] $G[S]$: the subgraph of $G$ induced by the vertices in $S$, where $S$
  is a subset of vertices of $G$.
  \item[3.] \textit{matching}, \textit{perfect matching}: A matching in a graph $G$ is a
  set of independent edges in $G$. A perfect matching $M$ in a graph
  $G$ is a matching such that every vertex of $G$ is incident to
  an edge of $M$.
  \item[4.] \textit{paired-dominating set}: A set $PD$ of vertices
  of $G$ is a paired-dominating set of $G$ if $PD$ is a dominating set of $G$ and if $G[PD]$
  contains at least one perfect matching.
  \item[5.] \textit{paired-domination number} $\gamma_{\rm p}(G)$: is the minimum cardinality
  of a paired-dominating set for a graph $G$.
  \item[6.] \textit{minimum paired-dominating set}: is a paired-dominating set of $G$
  with cardinality $\gamma_{\rm p}(G)$.
  \item[7.] $V(M)$: For a set $M$ of independent edges in a graph,
  $V(M)$ denotes the set of vertices being incident to edges of
  $M$.
  \item[8.] \textit{matched-paired-dominating set}: A set $MPD$ of independent edges in a graph
  $G$ is a matched-paired-dominating set of $G$ if $MPD$ is a perfect matching of
  $G[PD]$ induced by a paired-dominating set $PD$ of $G$.
  Note that $V(MPD)$ is a paired-dominating set $PD$ of $G$ and $MPD$ specifies a perfect matching
  of $G[PD]$.
  \item[9.] \textit{restricted vertex set} $\mathcal{R}$: The restricted vertex set $\mathcal{R}$
  is a subset of vertices in a graph and is a part of the input for the proposed problem in the
  paper. Any vertex in $\mathcal{R}$ is called \textit{restricted
  vertex} and the other is called \textit{free vertex}.
  \item[10.]\textit{maximum matched number} $\beta(G)$: For a matched-paired-dominating set
  $MPD$ of $G$, the matched number of $MPD$ is defined to be $|V(MPD)\cap \mathcal{R}|$.
  The maximum matched number $\beta(G)$ of $G$ is the largest matched number of a
  matched-paired-dominating set of $G$.
  \item[11.]\textit{maximum matched-paired-dominating set}: is a matched-paired-dominating set of a
  graph $G$ with matched number $\beta(G)$.
  \item[12.]\textit{paired-edge} $\langle u,v \rangle$: is an element in a matched-paired-dominating set $MPD$
  of a graph. We call $u$ the partner of $v$ in $MPD$.
  A paired-edge in $MPD$ is called \textit{full-paired-edge} if both of its vertices are
  in $\mathcal{R}$, is called \textit{semi-paired-edge} if its one vertex is
  in $\mathcal{R}$ but the other vertex is not in $\mathcal{R}$, and
  is called \textit{free-paired-edge} if both of its vertices are not in $\mathcal{R}$.
  \item[13.]\textit{canonical matched-paired-dominating set}: is a maximum matched-paired-dominating
  set of a graph $G$ with the least free-paired-edges.
  \item[14.]\textit{maximum matched-paired-domination problem}: Given a graph $G$ and a subset $\mathcal{R}$ of
  vertices in $G$, the problem is to find a canonical matched-paired-dominating set of
  $G$. Note that the proposed problem is a generalization of the paired-domination problem and it coincides with
  the classical paired-domination problem if $\mathcal{R}=\emptyset$.
  \item[15.]$(k,s,f)$-matched-paired-dominating set: is a matched-paired-dominating set $MPD$ of
  $G$ w.r.t. $\mathcal{R}$ satisfying that (1) $|MPD|=k+s+f$; (2) there are exactly $k$ full-paired-edges
  in $MPD$; (3) there are exactly $s$ semi-paired-edges in $MPD$,
  and (4) all other paired-edges in $MPD$ are free-paired-edges.
  \item[16.]$K_G(MPD)$, $S_G(MPD)$, $F_G(MPD)$: For a $(k,s,f)$-matched-paired-dominating set $MPD$ of a
  graph $G$, $K_G(MPD)$, $S_G(MPD)$, and $F_G(MPD)$ are defined to be
  the subsets of $MPD$ consisting of all full-paired-edges, all
  semi-paired-edges, and all free-paired-edges in $MPD$,
  respectively. That is, $|K_G(MPD)|=k$, $|S_G(MPD)|=s$, and $|F_G(MPD)|=f$.
  \item[17.]$G=G_L\oplus G_R$: $G=(V_L\cup V_R,E_L\cup E_R)$ is formed from $G_L=(V_L,E_L)$ and $G_R=(V_R,E_R)$ by a
  \textit{union operation}.
  \item[18.]$G=G_L\otimes G_R$: $G=(V,E)$ is formed from $G_L=(V_L,E_L)$ and
  $G_R=(V_R,E_R)$ by a \textit{joint operation}, where $V=V_L\cup V_R$
  and $E=E_L\cup E_R\cup \{uv|\forall u\in V_L$ and $\forall v\in V_R\}$.
  \item[19.]$I(G)$: is the set of isolated vertices in graph $G$.
  \item[20.]\textit{mixed paired-edges}: Let $G=G_L\otimes G_R$. A
  paired-edge in a matched-paired-dominating set of $G$ is called
  \textit{mixed} if one of its vertices is in $G_L$ and the other
  is in $G_R$.
  \item[21.] $\imath_L, \eta_R, f_R$: Let $G=G_L\otimes G_R$ with restricted vertex set
  $\mathcal{R}$, $\mathcal{R}_L=\mathcal{R}\cap V_L$, $\mathcal{R}_R=\mathcal{R}\cap
  V_R$, and let $CMPD_L$ be a canonical $(k_L,s_L,f_L)$-matched-paired-dominating
  set of $G_L-I(G_L)$. Define $\imath_L=|\mathcal{R}_L-V(CMPD_L)|$, $\eta_R=|\mathcal{R}_R|$,
  and $f_R=|V_R|-\eta_R$.
  \end{itemize}

\end{document}